\input amstex 
\documentstyle{amsppt}

\document
\magnification=1200
\NoBlackBoxes
\nologo
\vsize18cm
\centerline{\bf Quantized moduli spaces of the bundles on the elliptic curve}
\medskip 
\centerline{\bf and their applications.}
\bigskip
\centerline{\bf by A.V.Odesskii and B.L.Feigin}
\bigskip
\centerline{\bf Introduction}

\medskip

Let $\frak{b}$ be the Borel subalgebra of some Kac-Moody algebra $\frak{g}$
and $B$ the corresponding group. (If $\frak{b}$ is infinite-dimensional, then $B$ is a projective limit of finite-dimensional
groups). Let $\Cal M(\Cal E,\frak{b})$ be the moduli space of $B$-bundles on the elliptic curve $\Cal E$. There is the natural map $B\to H$, where $H$ is a Cartan subgroup, so each $B$-bundle defines an $H$-bundle. Note, that each $H$-bundle on $\Cal E$ has the characteristic class which lies in $H^2(\Cal E,L)\cong L$, where $L$ is the root lattice of Kac-Moody algebra $\frak{g}$. For each $n\in L$ we denote by $\Cal M_n(\Cal E,\frak{b})$ the corresponding connected component of $\Cal M(\Cal E,\frak{b})$.

In [3] we define the Hamiltonian structure on the manifold $\Cal M(\Cal E,\frak{b})$ and quantize the coordinate ring of the manifold $\Cal M_n(\Cal E,\frak{b})$ in the case when $n=\botshave{\sum_i}n_i\delta_i$, $\{\delta_i\}$ are simple positive roots, $\{n_i\}$ are positive and big enough. Actually, our construction works in case of arbitrary $\{n_i\}$. 

In this paper we consider the case when almost all $n_i$ are zero and some of them may be negative. 
Discuss first the case when all $\{n_i\}$ are zero. If $\frak{g}$ is finite-dimensional, then $\Cal M_0(\Cal E,\frak{b})$ is trivial. If $\frak{g}$ is an  affine Kac-Moody algebra, then $\Cal M_0(\Cal E,\frak{b})$ is infinite-dimensional with trivial Hamiltonian structure. Let $\frak{g}=\widehat{\frak{a}}$, where $\frak{a}$ is a  finite-dimensional semisimple Lie algebra and $P$ be a moduli space of $\frak{a}$-bundles on $\Cal E$. Then $\Cal M_0(\Cal E,\frak{b})$ is a space of maps $\varphi:D\to P$, where $D$ is a formal neighborhood of the origin in ${\Bbb {C}}$ such that $\varphi(0)$ is a trivial bundle. In {\bf \S5} we study the construction in this case.  In the case $\frak{g}=\widehat{sl_h}$ we write the quantization of $\Cal M_0(\Cal E,\frak{b})$ explicitly. It means that we construct an infinite set of commuting elements in some algebra. 

In {\bf \S4.2} we consider the case when $\frak{g}=sl_{h+1}$ and the set $(n_1,\dots,n_h)$ is $(m,0,0,\dots,0)$. The moduli space $\Cal M_{(m,0,\dots,0)}(\Cal E,\frak{b})$ has a projection on $\Cal E^h$ and a typical fiber is $(\Bbb {C}P^{m-1})^h$. After quantization we get a new construction of elliptic $R$-matrix.

In the {\bf \S4.4} we consider again the case $\frak{g}=sl_{h+1}$, when the set of $\{n_i\}$ is $(1,0,\dots,0,m)$. In this case our algebra contains the Sklyanin subalgebra $Q_{hm+1,h}(\Cal E,\tau)$ (see[2]).

In the {\bf \S4.5} again for $\frak{g}=sl_{h+1}$ we consider the set $(0,\dots,0,1,0,\dots,0)$. In this case we get the generalized elliptic $R$-matrix (compare [2], formula (6)). 

In the {\bf \S4.3} we consider the case $(-1,0,\dots,0,m)$ for $\frak{g}=sl_{h+1}$. It gives us the quantization of the coordinate ring of the Grassmannian of $h$-dimensional planes in $m-h+1$-dimensional space.

\medskip

Now we describe the contents of the paper.

 In {\bf \S1} we introduce the notations and construct the algebra $Q_{n,\Delta}(\Cal E,\tau)$. Here $\Delta$ is a root system of $\frak{g}$, $\tau\in \Bbb C$. For $\tau=0$ this algebra is the coordinate ring of the manifold $\Cal M_n(\Cal E,\frak{b})$ and for generic $\tau$ it is a quantization of this coordinate ring.

 In {\bf \S2} we construct some representations of the algebras $Q_{n,\Delta}(\Cal E,\tau)$ and in {\bf \S3} we study intertwining operators for these representations. Namely, we embed the algebra $Q_{n^\prime,\Delta}(\Cal E,-\tau)$ for some $n^\prime$ into the algebra of intertwining operators.

 In {\bf \S4} we apply the construction of the algebra $Q_{n,\Delta}(\Cal E,\tau)$ to produce some objects: the elliptic Belavin $R$-matrix, the quantization of the algebra of functions on the Grassmannian, the algebra $Q_{hm+1,h}(\Cal E,\tau)$ (see [2] for definition), some generalized elliptic $R$-matrix. 

In {\bf \S5} we consider the case of affine $\frak{g}$. We write down the explicit formula for commuting elements in the case $\frak{g}=\widehat{sl_h}$.

 In {\bf \S6} we construct the elliptic deformation of the Poison algebra and write down the explicit formula for commuting elements in this case (see also [5]). 

\newpage

\centerline{\bf \S1. Algebra $Q_{n,\Delta}(\Cal E,\tau)$.}
\medskip
{\bf 1.1. Notations.} Let $\Delta$ be the root system of the Kac-Moody algebra $\frak g$, $\Delta^+$ are the positive roots, $\{\delta_1,\dots,\delta_h\}$ are the simple positive roots (so $h$ is rank $\frak g$). We suppose that the set of simple positive roots in $\Delta$ is fixed. Let $(a_{i,j}), 1\leqslant i,j\leqslant h$ be the Cartan matrix, $a_{i,j}=\frac{2(\delta_i,\delta_j)}{(\delta_i,\delta_i)}$. Actually, for our construction of $Q_{n,\Delta}(\Cal E,\tau)$, we need only a linear space with a fixed basis $\{\delta_1,\dots,\delta_h\}$ and a scalar product $(\delta_i,\delta_j)$. But we don't know anything about the properties of the algebra $Q_{n,\Delta}(\Cal E,\tau)$ if $\{\delta_1,\dots,\delta_h\}$ is not a set of simple positive roots of the finite dimensional or affine Kac-Moody algebra. Let $L$ be a lattice with the basis $\{\delta_1,\dots,\delta_h\}$, $L\cong \Bbb Z^h$ and let $L^+=\{l_1\delta_1+\ldots+l_h\delta_h; l_1,\dots,l_h\geqslant 0\}\subset L$.

Let us $n: L\to\Bbb Z$ be homomorphism of the abelian groups, so $n_i=n(\delta_i)\in\Bbb Z$.

Let $\Cal E=\Bbb C\diagup\Gamma$ be an elliptic curve, where $\Gamma=\{m_1+m_2\eta; m_1,m_2\in \Bbb Z\}$ is a lattice, $\roman{Im}\,\eta>0$. For $m\in \Bbb Z$, $c\in\Bbb C$ we denote by $\xi_{m,c}$ a linear bundle on $\Cal E$ such that its sections are realized as functions on $\Bbb C$ with the following properties:$$f(z+1)=f(z); f(z+\eta)=e^{-2\pi i(mz+c)}f(z) \eqno(1)$$ It is clear that $\xi_{m,c}\cong\xi_{m^\prime,c^\prime}$ if and only if $m=m^\prime$ and $c-c^\prime\in\Gamma$. The space of holomorphic sections $H^0(\xi_{m,c})$ is realized as a space of holomorphic functions on $\Bbb C$ satisfying (1). We denote this space by $\Theta_{m,c}(\Gamma)$. It is easy to see that $\roman{dim}\,\Theta_{m,c}(\Gamma)=m$ if $m>0$, $\Theta_{m,c}(\Gamma)=0$ if $m<0$. If $m=0$ then $\Theta_{0,c}(\Gamma)=0$ for $c\notin\Gamma$. If $m=0$ and $c\in\Gamma$ then the bundle $\xi_{0,c}$ is trivial and $\roman{dim}\,H^0(\xi_{0,c})=1$. For $m>0$ the elements of $\Theta_{m,c}(\Gamma)$ are  called $\theta$-functions of order $m$. It is clear that for $m>0$ for different $c$ the spaces $\Theta_{m,c}(\Gamma)$ are isomorphic and can be identified by translation of $\Bbb C$. It is easy to check that every $\theta$-function of order $m$       has exactly $m$ zeros $\roman{mod}\,\Gamma$ and the sum of these zeros is equal to $c+\frac{1}{2}m$ $\roman{mod}\,\Gamma$. Let $\theta(z)=\botshave{\sum_{\alpha\in\Bbb Z}}(-1)^\alpha e^{2\pi i(\alpha z+\frac{\alpha(\alpha-1)}{2}\eta)}$. It is clear that $\theta(z)\in\Theta_{1,\frac{1}{2}}(\Gamma)$, $\theta(0)=0$, $\theta(-z)=-e^{-2\pi iz}\theta(z)$.

\medskip

{\bf 1.2. Algebra $F_\Delta(\Cal E,\tau)$.} Let $\tau\in\Bbb C$. We will construct an associative algebra $F_\Delta(\Cal E,\tau)$ that depends on 2 continues parameters: an  elliptic curve $\Cal E$ and a point $\tau\in\Bbb C$. The algebra $F_\Delta(\Cal E,\tau)$ is $L^+$-graded, so $F_\Delta(\Cal E,\tau)=\botshave{\bigoplus_{l\in L^+}}F_l$ and $F_l*F_{l^\prime}\subset F_{l+l^\prime}$. Here $*$ is the product in the algebra $F_\Delta(\Cal E,\tau)$. By definition, the space $F_l$ for $l=l_1\delta_1+\ldots+l_h\delta_h$ is a space of everywhere meromorphic functions $f(x_{1,1},\dots,x_{l_1,1};\dots;x_{1,h},\dots,x_{l_h,h};u_1,\dots,u_h)$ in $l_1+\ldots+l_h+h$ variables $\{x_{\alpha,i}, u_i; 1\leqslant i\leqslant h, 1\leqslant\alpha\leqslant l_i\}$. We assume, that these functions are symmetric with respect to each group of variables $\{x_{1,i},\dots,x_{l_i,i}\}$ for all $1\leqslant i\leqslant h$. Let $f\in F_l,g\in F_{l^\prime}$, where $l=l_1\delta_1+\ldots+l_h\delta_h, l^\prime=l_1^\prime\delta_1+\ldots+l_h^\prime\delta_h$. By definition, for the product $f*g\in F_{l+l^\prime}$ we have: $$f*g(x_{1,1},\dots,x_{l_1+l_1^\prime,1};\dots;x_{1,h},\dots,x_{l_h+l_h^\prime,h};u_1,\dots,u_h)=$$

$$\frac{1}{l_1!l_1^\prime!\ldots l_h!l_h^\prime!}\sum\Sb\sigma_1\in S_{l_1+l_1^\prime}\\ \\\dots\dots\dots\dots\dots\\ \\\sigma_h\in S_{l_h+l_h^\prime}\endSb f(x_{\sigma_1(1),1},\dots,x_{\sigma_1(l_1),1};\dots;x_{\sigma_h(1),h},\dots,x_{\sigma_h(l_h),h};u_1,\dots,u_h)\times$$

$$g(x_{\sigma_1(l_1+1),1},\dots,x_{\sigma_1(l_1+l_1^\prime),1};\dots;x_{\sigma_h(l_h+1),h},\dots,x_{\sigma_h(l_h+l_h^\prime),h};u_1^{(\tau)},\dots,u_h^{(\tau)})\times \eqno(2)$$

$$\botshave{\prod\Sb1\leqslant i,j\leqslant h,\\ 1\leqslant\alpha\leqslant l_i,\\ l_j+1\leqslant\beta\leqslant l_j+l_j^\prime\endSb}\frac{\theta(x_{\sigma_i(\alpha),i}-x_{\sigma_j(\beta),j}-(\delta_i,\delta_j)\tau)}{\theta(x_{\sigma_i(\alpha),i}-x_{\sigma_j(\beta),j})}$$
Here $u_i^{(\tau)}=u_i-2(l,\delta_i)\tau$ for $1\leqslant i\leqslant h$.

{\bf Proposition 1.} {\it The formula (2) defines in the space $F_\Delta(\Cal E,\tau)=\botshave{\bigoplus_{l\in L^+}}F_l$ the structure of the associative algebra.}

{\bf Proof.} By direct calculation.

{\bf Remarks. 1.} Let $\tau=0$. In this case the algebra $F_\Delta(\Cal E,0)$ is commutative and does not depend on $\Cal E$. The formula for $f*g$ takes such a form:
$$f*g(x_{1,1},\dots,x_{l_1+l_1^\prime,1};\dots;x_{1,h},\dots,x_{l_h+l_h^\prime,h};u_1,\dots,u_h)=$$

$$\frac{1}{l_1!l_1^\prime!\ldots l_h!l_h^\prime!}\sum\Sb\sigma_1\in S_{l_1+l_1^\prime}\\ \\\dots\dots\dots\dots\dots\\ \\\sigma_h\in S_{l_h+l_h^\prime}\endSb f(x_{\sigma_1(1),1},\dots,x_{\sigma_1(l_1),1};\dots;x_{\sigma_h(1),h},\dots,x_{\sigma_h(l_h),h};u_1,\dots,u_h)\times$$

$$g(x_{\sigma_1(l_1+1),1},\dots,x_{\sigma_1(l_1+l_1^\prime),1};\dots;x_{\sigma_h(l_h+1),h},\dots,x_{\sigma_h(l_h+l_h^\prime),h};u_1,\dots,u_h) \eqno(3)$$

Let $P_l\subset F_l$ be a space of functions which does not depend on $\{u_1,\dots,u_h\}$ and are polynomials in variables $\{x_{\alpha,i}; 1\leqslant i\leqslant h, 1\leqslant\alpha\leqslant l_i\}$. It is clear that $P=\botshave{\bigoplus_{l\in L^+}}P_l$ is subalgebra in $F_\Delta(\Cal E,0)$. We have: $P_l=S^{l_1}P_{\delta_1}\otimes S^{l_2}P_{\delta_2}\otimes\ldots\otimes S^{l_h}P_{\delta_h}$, so $P=S^*P_{\delta_1}\otimes S^*P_{\delta_2}\otimes\ldots\otimes S^*P_{\delta_h}$. Here $P_{\delta_i}$ is a space of polynomials in one variable. The formula (3) defines the usual product in this tensor product of symmetric algebras.

{\bf 2.} By definition $F_0$ is a space of all meromorphic functions in variables $\{u_1,\dots,u_h\}$. For $f,g\in F_0$ we have (see (2)) $f*g=fg$. It is the usual product of functions. So $F_0$ is a field. We have $F_0*F_l\subset F_l$ and $F_l*F_0\subset F_l$ so we have two structures of a $F_0$-vector space on $F_l$. These two structures are connected by the following: for $f\in F_0$ and $g\in F_l$ we have: $g*f(u_1,\dots,u_h)=f(u_1-2(l,\delta_1)\tau,\dots,u_h-2(l,\delta_h)\tau)*g$ (see (2)). For $\tau=0$, the algebra $F_\Delta(\Cal E,0)$ is a $F_0$-algebra.

\medskip

{\bf 1.3. Algebra $Q_{n,\Delta}(\Cal E,\tau)$.} For every linear form $n: L\to \Bbb Z$ we  will construct a subalgebra $Q_{n,\Delta}(\Cal E,\tau)$ in the algebra $F_\Delta(\Cal E,\tau)$. By definition $Q_{n,\Delta}(\Cal E,\tau)=\botshave{\bigoplus_{l\in L^+}}Q_l$, here $Q_l\subset F_l$ is a space of functions satisfied following conditions:

1. $f(x_{1,1},\dots,x_{l_h,h};u_1,\dots,u_h)$ as a function in variables $\{x_{1,1},\dots,x_{l_h,h}\}$ is holomorphic outside the divisors $\{x_{\alpha,i}-x_{\beta,j}=0; i\ne j, (\delta_i,\delta_j)\ne0\}$ and has a pole of order $\leqslant1$ on these divisors.

2. For each $x_{\alpha,i}$, $1\leqslant i\leqslant h$, $1\leqslant\alpha\leqslant l_i$  $f$ satisfies (1) for $m=n_i$, $c=u_i-(\delta_i,l)\tau$ as a function in $x_{\alpha,i}$. Here $n_i=n(\delta_i)$. So we have: $$f(x_{1,1},\dots,x_{\alpha,i}+1,\dots,x_{l_h,h};u_1,\dots,u_h)=f(x_{1,1},\dots,x_{l_h,h};u_1,\dots,u_h);$$

$$f(x_{1,1},\dots,x_{\alpha,i}+\eta,\dots,x_{l_h,h};u_1,\dots,u_h)=$$
$$e^{-2\pi i(n_ix_{\alpha,i}+u_i-(\delta_i,l)\tau)}f(x_{1,1},\dots,x_{l_h,h};u_1,\dots,u_h)$$

3. Let $1\leqslant i,j\leqslant h$ such that $a_{i,j}<0$. Let $1\leqslant\alpha_1,\dots,\alpha_{-a_{i,j}+1}\leqslant l_i$, $1\leqslant\beta\leqslant l_j$. Then $f(x_{1,1},\dots,u_h)=0$ on the affine subspace of codimension $-a_{i,j}+1$  defined by the following relations: $$x_{\alpha_1,i}-x_{\alpha_2,i}=(\delta_i,\delta_i)\tau,
x_{\alpha_2,i}-x_{\alpha_3,i}=(\delta_i,\delta_i)\tau,\dots,
x_{\alpha_{-a_{i,j}},i}-x_{\alpha_{-a_{i,j}+1},i}=(\delta_i,\delta_i)\tau,$$

$$x_{\alpha_{-a_{i,j}+1},i}-x_{\beta,j}=(\delta_i,\delta_j)\tau,
x_{\beta,j}-x_{\alpha_1,i}=(\delta_i,\delta_j)\tau.$$

4. Let $\{\delta_{i_1},\dots,\delta_{i_p}\}$ generates an irreducible component of the root system $\Delta$; $\mu_1,\dots,\mu_h;\lambda_1,\dots,\lambda_h;\nu\in\Bbb C$ and $\mu_j=\lambda_j=0$ if $j\notin\{i_1,\dots,i_p\}$, $\mu_{i_\alpha}=-\nu$, $\lambda_{i_\alpha}=n_{i_\alpha}\nu$ for $1\leqslant\alpha\leqslant p$. Then $f(x_{1,1}+\mu_1,\dots,x_{l_1,1}+\mu_1;\dots;x_{1,h}+\mu_h,\dots,x_{l_h,h}+\mu_h;u_1+\lambda_1,\dots,u_h+\lambda_h)=f(x_{1,1},\dots,u_h)$.

{\bf Proposition 2.} {\it The space $Q_{n,\Delta}(\Cal E,\tau)$ is closed with respect to the product $*$ defined by (2). So $Q_{n,\Delta}(\Cal E,\tau)$ is a subalgebra of the algebra $F_\Delta(\Cal E,\tau)$.}

{\bf Proof.} By direct calculation it is easy to see that the product $*$ preserves each of the properties 1-4.

{\bf Proposition 3.} {\it Let $\Delta_1\subset\Delta$ is a root subsystem generated by some subset of $\{\delta_1,\dots,\delta_h\}$. Then $F_{\Delta_1}(\Cal E,\tau)\subset F_\Delta(\Cal E,\tau)$, $Q_{m,\Delta_1}(\Cal E,\tau)\subset Q_{n,\Delta}(\Cal E,\tau)$ as graded subalgebras. Here $m$ is a restriction of $n$ to the sublattice generated by $\Delta_1$.}

{\bf Proof.} It is evident from the definitions.

{\bf Remarks. 1.} $Q_0$ is a field of meromorphic functions $f(u_1,\dots,u_h)$ satisfying the property 4.

{\bf 2.} We have $Q_0*Q_l\subset Q_l$ and $Q_l*Q_0\subset Q_l$ so $Q_l$ has two structures of a $Q_0$-vector space. These two structures are connected in the following way: for $f\in Q_0$ and $g\in Q_l$ we have; $g*f(u_1,\dots,u_h)=f(u_1-2(l,\delta_1)\tau,\dots,u_h-2(l,\delta_h)\tau)*g$ (see (2)), so the dimension of these spaces is the same.

{\bf 3.} From the properties 1 and 2 of the elements from the space $Q_l$ it follows that $\roman{dim}_{Q_0}Q_l$ is finite for each $l\in L^+$.

{\bf 4.} For $\tau=0$ the algebra $Q_{n,\Delta}(\Cal E,\tau)$ is a commutative $Q_0$-algebra.

{\bf 5.} $Q_{\delta_i}$ consists of the functions $f(x;u_1,\dots,u_h)$ which are  holomorphic with respect to $x$ and satisfy the properties 2 and 4. So $\roman{dim}_{Q_0}Q_{\delta_i}=n_i$ if $n_i>0$ and $Q_{\delta_i}=0$ if $n_i\leqslant0$.

{\bf 6.} $Q_{\alpha\delta_i}$ consists of the functions $f(x_1,\dots,x_\alpha;u_1,\dots,u_h)$ which are holomorphic and symmetric with respect to $\{x_1,\dots,x_\alpha\}$ and satisfy the properties 2 and 4. It is clear that $\roman{dim}_{Q_0}Q_{\alpha\delta_i}=\frac{n_i(n_i+1)\ldots(n_i+\alpha-1)}{\alpha!}$ if $n_i>0$ and $Q_{\alpha\delta_i}=0$ if $n_i\leqslant0$.

{\bf 7.} If $(\delta_i,\delta_j)=0$ then $Q_{\alpha\delta_i+\beta\delta_j}\cong Q_{\alpha\delta_i}\otimes_{Q_0}Q_{\beta\delta_j}$.

{\bf 8.} If $(\delta_i,\delta_j)\ne0$ (and $i\ne j$), than $Q_{\delta_i+\delta_j}$ becomes bigger then $Q_{\delta_i}\otimes_{Q_0}Q_{\delta_j}$ because of the pole part of the functions from $Q_{\delta_i+\delta_j}$. So we have: $\roman{dim}_{Q_0}Q_{\delta_i+\delta_j}=n_in_j+n_i+n_j$, if $n_i, n_j\geqslant0$. If $n_i$ or $n_j<0$, then $Q_{\delta_i+\delta_j}=0$.

\newpage

\centerline{\bf \S2. Representations of the algebras $Q_{n,\Delta}(\Cal E,\tau)$.}
\medskip
Let $p\in L^+$, $p=p_1\delta_1+\ldots+p_h\delta_h$. Let $A_{p,\Delta}(\Cal E,\tau)$ be the algebra generated by $\{e_{\alpha,i}, e_{\alpha,i}^{-1}; 1\leqslant i\leqslant h, 1\leqslant\alpha\leqslant p_i\}$ and $\{\varphi(y_{1,1},\dots,y_{p_1,1};\dots;y_{1,h},\dots,y_{p_h,h};u_1,\dots,u_h)\}$, where $\varphi$ is any meromorphic function in variables $\{y_{\alpha,i}, u_i; 1\leqslant i\leqslant h, 1\leqslant\alpha\leqslant p_i\}$ such that as a function of $\{y_{\alpha,i}\}$ it has poles only on the divisors of the form $\{y_{\alpha,i}-y_{\beta,j}-(\delta_i,\delta_j)\tau+\mu(\delta_i,\delta_i)\tau-\nu(\delta_j,\delta_j)\tau=0; \mu, \nu\in\Bbb Z, (\delta_i,\delta_j)\ne0\}$. We assume such the following relations hold: $$e_{\alpha,i}e_{\beta,j}=-e^{2\pi i(y_{\beta,j}-y_{\alpha,i})}\frac{\theta(y_{\alpha,i}-y_{\beta,j}-(\delta_i,\delta_j)\tau)}{\theta(y_{\beta,j}-y_{\alpha,i}-(\delta_i,\delta_j)\tau)}e_{\beta,j}e_{\alpha,i}\eqno(4)$$
$e_{\alpha,i}y_{\beta,j}=y_{\beta,j}e_{\alpha,i}$, here $i\ne j$ or $\alpha\ne\beta$;
$$e_{\alpha,i}e_{\alpha,i}^{-1}=1, e_{\alpha,i}y_{\alpha,i}=(y_{\alpha,i}+(\delta_i,\delta_i)\tau)e_{\alpha,i}, e_{\alpha,i}u_j=(u_j-2(\delta_i,\delta_j)\tau)e_{\alpha,i};$$
$$[y_{\alpha,i},y_{\beta,j}]=[y_{\alpha,i},u_j]=[u_i,u_j]=0.$$
Let us define the map $x: Q_{n,\Delta}(\Cal E,\tau)\to A_{p,\Delta}(\Cal E,\tau)$ by the following. For $f\in Q_0$, we put $x(f)=f$. Recall that $Q_0\subset A_{p,\Delta}(\Cal E,\tau)$ by definition. For $f(x;u_1,\dots,u_h)\in Q_{\delta_i}$ (if $n_i>0$ ), we put $x(f)=\botshave{\sum_{1\leqslant\alpha\leqslant p_i}}f(y_{\alpha,i};u_1,\dots,u_h)e_{\alpha,i}$. In general case for $f(x_{1,1},\dots,x_{l_h,h};u_1,\dots,u_h)\in Q_l$, $l=l_1\delta_1+\ldots+l_h\delta_h$, we define: $$x(f)=\sum\Sb\varphi_{1,1},\dots,\varphi_{p_h,h}\geqslant0,\\ \varphi_{1,i}+\ldots+\varphi_{p_i,i}=l_i\text{ for each }1\leqslant i\leqslant h\endSb B_{\varphi_{1,1},\dots,\varphi_{p_h,h}}e_{1,1}^{\varphi_{1,1}}\ldots e_{p_1,1}^{\varphi_{p_1,1}}\ldots e_{1,h}^{\varphi_{1,h}}\ldots e_{p_h,h}^{\varphi_{p_h,h}} \eqno(5)$$ Here $B_{\varphi_{1,1},\dots,\varphi_{p_h,h}}$ is a function in variables $\{y_{\alpha,i},u_i\}$, $$B_{\varphi_{1,1},\dots,\varphi_{p_h,h}}=$$
 $$f(y_{1,1},y_{1,1}+(\delta_1,\delta_1)\tau,\dots,y_{1,1}+(\varphi_{1,1}-1)(\delta_1,\delta_1)\tau,\dots,$$
$$y_{p_h,h},y_{p_h,h}+(\delta_h,\delta_h)\tau,\dots,y_{p_h,h}+(\varphi_{p_h,h}-1)(\delta_h,\delta_h)\tau;u_1,\dots,u_h)\times$$
$$\prod\Sb1\leqslant i,j\leqslant h,\\ 1\leqslant\alpha\leqslant p_i,1\leqslant\beta\leqslant p_j,\\ 0\leqslant\mu\leqslant\varphi_{\alpha,i}-1, 0\leqslant\nu\leqslant\varphi_{\beta,j}-1;\\ i\leqslant j;\text{ if $i=j$, then $\alpha\leqslant\beta$;}\\\text{ if $i=j,\alpha=\beta$, then $\mu<\nu$}\endSb\frac{\theta(y_{\alpha,i}+\mu(\delta_i,\delta_i)\tau-y_{\beta,j}-\nu(\delta_j,\delta_j)\tau)}{\theta(y_{\alpha,i}+\mu(\delta_i,\delta_i)\tau-y_{\beta,j}-\nu(\delta_j,\delta_j)\tau-(\delta_i,\delta_j)\tau)}$$

{\bf Proposition 4.} {\it $x$ is homomorphism of algebras.}

{\bf Proof.} It is direct checking using the formulas (2) and (4).

{\bf Remarks. 1.} If $\varphi_{\alpha,i}=0$ for some $\alpha$ and $i$, then $f$ does not contain variable $y_{\alpha,i}$ in corresponding part of the formula (5). For $\varphi_{\alpha,i}>0$, $f$ contains variables $y_{\alpha,i}+\mu(\delta_i,\delta_i)\tau$ for $0\leqslant\mu\leqslant\varphi_{\alpha,i}-1$.

{\bf 2.} It is possible to construct representations of the algebra $Q_{n,\Delta}(\Cal E,\tau)$ using homomorphism $x$. For this let us construct the representation of the algebra $A_{p,\Delta}(\Cal E,\tau)$ with diagonal action of the elements $\{y_{\alpha,i},u_i\}$. The basis of this representation consists of  monomials of the elements $\{e_{\alpha,i}\}$. The homomorphism $x$ defines the representation of the algebra $Q_{n,\Delta}(\Cal E,\tau)$ in this space.

\newpage

\centerline{\bf \S3. Intertwining operators.}
\medskip

In this paragraph we assume that $\Delta$ is a root system of the finite dimension Lie algebra. The intertwining operators for the representations of the algebra $Q_{n,\Delta}(\Cal E,\tau)$ are constructed in the following way (see also [6], \S3). Let $\sigma$ be an automorphism of the algebra $A_{p,\Delta}(\Cal E,\tau)$. We say that an element $C\in A_{p,\Delta}(\Cal E,\tau)$ $\sigma$-commutes with the image $x(Q_{n,\Delta}(\Cal E,\tau))$, if for each element $f\in Q_{n,\Delta}(\Cal E,\tau)$ we have: $Cx(f)=(\sigma x(f))C$ in the algebra $A_{p,\Delta}(\Cal E,\tau)$. Let $\pi$ be some representation of the algebra $Q_{n,\Delta}(\Cal E,\tau)$ constructed by the homomorphism $x$ (see remark 2 in \S2). It is clear that $C$ defines the intertwining operator between the representations $\pi$ and $\pi^\sigma$ of the algebra $Q_{n,\Delta}(\Cal E,\tau)$. Let $C_1\in A_{p,\Delta}(\Cal E,\tau)$ $\sigma_1$-commutes and $C_2\in A_{p,\Delta}(\Cal E,\tau)$ $\sigma_2$-commutes with $x(Q_{n,\Delta}(\Cal E,\tau))$. It is clear that $(\sigma_2C_1)C_2$ $\sigma_2\sigma_1$-commutes with $x(Q_{n,\Delta}(\Cal E,\tau))$. This product coincides with the product of the intertwining operators.

Let $L_{\Bbb C}=L\otimes_{\Bbb Z}{\Bbb C}$ be a $\Bbb C$-vector space with a basis $\{\delta_1,\dots,\delta_h\}$. Let $\mu\in L_{\Bbb C}$, $\mu=\mu_1\delta_1+\ldots+\mu_h\delta_h$. Let us define an automorphism $T_{\mu}$ of the algebra $A_{p,\Delta}(\Cal E,\tau)$ by the following: $T_{\mu}(e_{\alpha,i})=e_{\alpha,i}$, $T_{\mu}(y_{\alpha,i})=y_{\alpha,i}$, $T_{\mu}(u_i)=u_i+\mu_i$. We will look for the elements of the algebra $A_{p,\Delta}(\Cal E,\tau)$ that $T_{\mu}$-commutes with $x(Q_{n,\Delta}(\Cal E,\tau))$. The algebra $A_{p,\Delta}(\Cal E,\tau)$ is $L$-graded, if we assume $\roman{deg}\,y_{\alpha,i}=\roman{deg}\,u_i=0, \roman{deg}\,e_{\alpha,i}=\delta_i$ for each $1\leqslant i\leqslant h$, $1\leqslant\alpha\leqslant p_i$. It is clear that $A_{p,\Delta}(\Cal E,\tau)=\botshave{\bigoplus_{l\in L}}A_l$, here $A_l$ is the space of elements of degree $l$ and $A_lA_{l^\prime}\subset A_{l+l^\prime}$ in the algebra $A_{p,\Delta}(\Cal E,\tau)$. It is clear that the homomorphism $x$ preserves the grading, so $x(Q_l)\subset A_l$. Let $\omega: L\to L_{\Bbb C}$ be a homomorphism of abelian groups. We denote by $A_{p,\Delta}^{\omega}(\Cal E,\tau)$ the associative algebra such that as a linear space it is isomorphic to $A_{p,\Delta}(\Cal E,\tau)$ and the product $\circ$ is defined by the following: $C_1\circ C_2=(T_{\omega(l_2)}C_1)C_2$. Here $C_1\in A_{l_1}, C_2\in A_{l_2}$ and $\circ$ is the product in the algebra $A_{p,\Delta}^\omega(\Cal E,\tau)$. We will construct a  homomorphism $$y: Q_{n^\prime,\Delta}(\Cal E,-\tau)\to A^\omega_{p,\Delta}(\Cal E,\tau)$$ for some $n^\prime$ and $\omega$ such that for each elements $f\in Q_{n,\Delta}(\Cal E,\tau)$,  $g\in Q_{n^\prime,\Delta}(\Cal E,-\tau)$ we will have: $$y(g)x(f)=(T_{-\omega(l)}x(f))y(g)$$ in the algebra $A_{p,\Delta}(\Cal E,\tau)$. Here $g\in Q^\prime_l$, $Q_{n^\prime,\Delta}(\Cal E,-\tau)=\botshave{\bigoplus_{l\in L^+}}Q^\prime_l$ is the decomposition of $Q_{n^\prime,\Delta}(\Cal E,-\tau)$ for $L^+$-grading. We will have $y(Q^\prime_l)\subset A_{-l}$.

Let $n^\prime: L\to \Bbb Z$ be such homomorphism that $$n^\prime(\delta_i)=\botshave{\sum_{1\leqslant j\leqslant h}}a_{j,i}p_j-n(\delta_i)$$ We denote $n^\prime_i=n^\prime(\delta_i)$ for all $1\leqslant i\leqslant h$.

Let $A^\prime_{p,\Delta}(\Cal E,-\tau)$ be the algebra with generators $\{e^\prime_{\alpha,i}, e^{\prime-1}_{\alpha,i}; 1\leqslant i\leqslant h, 1\leqslant\alpha\leqslant p_i\}$ and $\{\varphi(y^\prime_{1,1},\dots,y^\prime_{p_1,1};\dots;y^\prime_{1,h},\dots,y^\prime_{p_h,h};u^\prime_1,\dots,u^\prime_h)\}$, where $\varphi$ are meromorphic functions. We assume the the relations are as in (4), with $\tau$ replaced by $-\tau$ and $e_{\alpha,i}, y_{\alpha,i}, u_i$ replaced by $e^\prime_{\alpha,i}, y^\prime_{\alpha,i}, u^\prime_i$ respectively. It is clear, that there is a homomorphism  $x^\prime: Q_{n^\prime,\Delta}(\Cal E,-\tau)\to A^\prime_{p,\Delta}(\Cal E,-\tau)$ defined by (5), where $Q_{n,\Delta}(\Cal E,\tau)$ is replaced by $Q_{n^\prime,\Delta}(\Cal E,-\tau)$ and $A_{p,\Delta}(\Cal E,\tau)$ is replaced by $A^\prime_{p,\Delta}(\Cal E,-\tau)$. Here $Q^\prime_0$ consists of meromorphic functions in variables $\{u^\prime_1,\dots,u^\prime_h\}$.

Let $\omega: L\to L_{\Bbb C}$ be such homomorphism that $$\omega(\delta_i)=-\sum_{1\leqslant j\leqslant h}2(\delta_i,\delta_j)\tau\cdot\delta_j$$

{\bf Proposition 5.} {\it Following formulas define the homomorphism of the algebras $\varkappa: A^\prime_{p,\Delta}(\Cal E,-\tau)\to A^\omega_{p,\Delta}(\Cal E,\tau)$:
$$\varkappa(y^\prime_{\alpha,i})=y_{\alpha,i}$$
$$\varkappa(u^\prime_i)=-u_i-\sum\Sb1\leqslant j\leqslant h,\\ 1\leqslant\beta\leqslant p_j\endSb a_{j,i}y_{\beta,j} \eqno(6)$$
$$\varkappa(e^\prime_{\alpha,i})=\prod\Sb j\ne i, 1\leqslant j\leqslant h,\\ 1\leqslant\beta\leqslant p_j, 0\leqslant\delta\leqslant-a_{j,i}-1\endSb(\theta(y_{\alpha,i}-y_{\beta,j}-((\delta_i,\delta_j)+(\delta_i,\delta_i)+\delta(\delta_j,\delta_j))\tau)e^{\pi iy_{\beta,j}})\times$$
$$\prod\Sb1\leqslant\beta\leqslant p_i\\ \beta\ne\alpha\endSb\frac{e^{-2\pi iy_{\beta,i}}}{\theta(y_{\alpha,i}-y_{\beta,i}-(\delta_i,\delta_i)\tau)\theta(y_{\alpha,i}-y_{\beta,i})}\cdot e^{-1}_{\alpha,i}$$}

{\bf Proof.} It is a direct checking of the relations (4).

Let $T^\prime_\mu$ be an automorphism of the algebra $A^\prime_{p,\Delta}(\Cal E,-\tau)$ such that $T^\prime_\mu e^\prime_{\alpha,i}=e^\prime_{\alpha,i}$, $T^\prime_\mu y^\prime_{\alpha,i}=y^\prime_{\alpha,i}$, $T^\prime_\mu u^\prime_i=u^\prime_i+\mu_i$. Here $\mu\in L_{\Bbb C}$, $\mu=\mu_1\delta_1+\ldots+\mu_h\delta_h$. We put: $$\mu_i=(-(n^\prime_i+3)(\delta_i,\delta_i)+(\delta_i,p))\tau+\frac{1}{2}\sum_{1\leqslant j\leqslant h}a_{j,i}p_j$$

{\bf Proposition 6.} {\it Let $y: Q_{n^\prime,\Delta}(\Cal E,-\tau)\to A^\omega_{p,\Delta}(\Cal E,\tau)$ be the composition of the homomorphisms $x^\prime$, $T^\prime_\mu$ and $\varkappa$, so $y(g)=\varkappa T^\prime_\mu x^\prime(g)$ for $g\in Q_{n^\prime,\Delta}(\Cal E,-\tau)$. Then for each $g\in Q_{n^\prime,\Delta}(\Cal E,-\tau)$, $f\in Q_{n,\Delta}(\Cal E,\tau)$ there is a relation $y(g)x(f)=(T_{-\omega(l)}x(f))y(g)$ in the algebra $A_{p,\Delta}(\Cal E,\tau)$. Here $g\in Q^\prime_l$.}

{\bf Proof.} See [6], \S3 for the case of strict dominant $n$ and $n^\prime$ (this means $n(\delta_i)>0$ and $n^\prime(\delta_i)>0$). In general case the proof is similar.

{\bf Remarks. 1.} This construction for the case of strict dominant $n$ and $n^\prime$ was done in [6], \S3, propositions 3 and 4. We have extended this construction for general case above.

{\bf 2.} In [6], \S3, Proposition 5 we have constructed another elements $K_d\in A_{p,\Delta}(\Cal E,\tau)$ that $T_\nu$-commute with $x(Q_{n,\Delta}(\Cal E,\tau))$. Here $d$ is such that $\frac{d}{(\delta_i,\delta_i)}\in\Bbb N$ for each $1\leqslant i\leqslant h$, $\nu=-d\tau(n^\prime_1\delta_1+\ldots+n^\prime_h\delta_h)$. In [6] we considered only the case of strict dominant $n$ and $n^\prime$, but it is clear that the Proposition 5 holds for general case.

{\bf 3.} In [6] we used slightly different notations (compare [6], \S1 and \S1 of this paper). In notations of [6] we had to divide by $n_i$ (see for example [6], formula (1)) but it is not good if we want to consider the case $n_i=0$ for some $i$. So we have to change the notations.

\newpage

\centerline{\bf \S4. Examples and applications.}

\medskip

{\bf 4.1. Case of finite root system $\Delta$ and dominant $n$.}

This case was studied in [6]. We remind the main result about the structure of the algebra $Q_{n,\Delta}(\Cal E,\tau)$. 

Let $\frak g=\frak g^-\oplus\frak h\oplus\frak g^+$ be the Cartan decomposition, $\{g_\gamma; \gamma\in\Delta^+\}$ be a basis in $\frak g^+$, $\varphi_{\gamma_1,\gamma_2}\in \Bbb C$ be the structure constants, so $[g_{\gamma_1},g_{\gamma_2}]=\varphi_{\gamma_1,\gamma_2}g_{\gamma_1+\gamma_2}$. We denote by $\widehat{\frak g^+}(\Cal E,n)$ the following $L^+$-graded Lie algebra over field $Q_0$: $$\widehat{\frak g^+}(\Cal E,n)=\botshave{\bigoplus_{\gamma\in\Delta^+}}\Theta_{n(\gamma),u(\gamma)}(\Gamma)$$ as a linear space. Here $u(\gamma)=u_1\gamma_1+\ldots+u_h\gamma_h$ for $\gamma=\gamma_1\delta_1+\ldots+\gamma_h\delta_h$. The commutator for $\mu\in\Theta_{n(\gamma),u(\gamma)}(\Gamma)$, $\lambda\in\Theta_{n(\gamma^\prime),u(\gamma^\prime)}(\Gamma)$ is by definition $[\mu,\lambda]=\varphi_{\gamma,\gamma^\prime}\mu\cdot\lambda$. Here $\mu\cdot\lambda$ is the usual product of functions. We remind that by definition $\Theta_{n(\gamma),u(\gamma)}(\Gamma)=0$ if $n(\gamma)=0$ (see \S1.1). It is clear that  $\widehat{\frak g^+}(\Cal E,n)$ is a subalgebra of the current algebra of $\frak g^+$ over field $Q_0$.

{\bf Proposition 7.} {\it The algebra $Q_{n,\Delta}(\Cal E,\tau)$, defined for finite root system $\Delta$ and dominant $n$ is a deformation of the universal enveloping algebra $U(\widehat{\frak g^+}(\Cal E,n))$ in the class of $L^+$-graded associative algebras.}

{\bf Remark.} From this proposition it follows that the Hilbert function of the algebra $Q_{n,\Delta}(\Cal E,\tau)$ for finite $\Delta$ and dominant $n$ is $$\botshave{\sum_{l\in L^+}}\roman{dim}_{Q_0}Q_l w^l=\botshave{\prod_{\delta\in\Delta^+}}(1-w^\delta)^{-n(\delta)}$$ Here $w^l=w_1^{l_1}\ldots w_h^{l_h}$ for $l=l_1\delta_1+\ldots+l_h\delta_h$; $w_1,\dots,w_h$ are formal parameters.

\medskip

{\bf 4.2. The Zamolodchikove algebra for the elliptic Belavin $R$-matrix.}

Let $\Delta$ be root system $A_h$, so $(\delta_i,\delta_i)=2$, $1\leqslant i\leqslant h$; $(\delta_i,\delta_{i+1})=-1$, $1\leqslant i<h$; $(\delta_i,\delta_j)=0$ if $|i-j|>1$. Let $n(\delta_1)=m>0$, $n(\delta_2)=\ldots=n(\delta_h)=0$. The Hilbert function of the algebra $Q_{n,\Delta}(\Cal E,\tau)$ has a form (see \S4.1): $\botshave{\prod_{\gamma\in\Delta^+}}(1-w^\gamma)^{-n(\gamma)}=(1-w_1)^{-m}(1-w_1w_2)^{-m}\ldots(1-w_1\ldots w_h)^{-m}$. So the algebra $Q_{n,\Delta}(\Cal E,\tau)$ is a deformation of the polynomial ring in $hm$ variables: $m$ variables of degree $\delta_1+\ldots+\delta_i$ for each $1\leqslant i\leqslant h$. It is easy to check that the space $Q_{\delta_1+\ldots+\delta_i}$ consists of functions of the form:
$$f(x_{1,1};\dots;x_{1,i};u_1,\dots,u_h)=\frac{\botshave{\prod_{1\leqslant\alpha<i}}\theta(x_{1,\alpha}-x_{1,\alpha+1}-\botshave{\sum_{\alpha<\nu\leqslant i}}u_{\nu}+\tau)}{\botshave{\prod_{1\leqslant\alpha<i}}\theta(x_{1,\alpha}-x_{1,\alpha+1})}\varphi(x_{1,1}),$$
where $\varphi(x)\in\Theta_{m,u_1+\ldots+u_i}(\Gamma)$.

{\bf Proposition 8.} {\it The space $Q_{\delta_1+\ldots+\delta_i}$ generates the algebra $Q_m(\Cal E,\frac{2}{m}\tau)$ for each $1\leqslant i\leqslant h$. The commutation relations between $Q_{\delta_1+\ldots+\delta_i}$ and $Q_{\delta_1+\ldots+\delta_j}$, $i\ne j$, $1\leqslant i,j\leqslant h$ are defined by the elliptic $R$-matrix $R_m(\Cal E,\frac{2}{m}\tau)$ (for the definitions of $Q_m(\Cal E,\tau)$ and $R_m(\Cal E,\tau)$ see [1]).}

{\bf Proof.} It is enough to proof the proposition for the image of the homomorphism $x: Q_{n,\Delta}(\Cal E,\tau)\to A_{p,\Delta}(\Cal E,\tau)$, because for large $\{p_i\}$ this homomorphism is injection. Let us denote by $z_{\varphi,i}$ the image of the element from the space $Q_{\delta_1+\ldots+\delta_i}$ for a function  $\varphi(x)\in\Theta_{m,u_1+\ldots+u_i}(\Gamma)$. We denote by $v_i=u_1+\ldots+u_i$ for $1\leqslant i\leqslant h$. The formula (5) gives us the following expression for $z_{\varphi,i}\in A_{p,\Delta}(\Cal E,\tau)$:
$$z_{\varphi,i}=\sum\Sb 1\leqslant\alpha_1\leqslant p_1\\ \\\dots\dots\dots\dots\\ \\1\leqslant\alpha_i\leqslant p_i\endSb\frac{\botshave{\prod_{1\leqslant\nu<i}}\theta(y_{\alpha_\nu,\nu}-y_{\alpha_{\nu+1},\nu+1}+v_\nu-v_i+\tau)}{\botshave{\prod_{1\leqslant\nu<i}}\theta(y_{\alpha_\nu,\nu}-y_{\alpha_{\nu+1},\nu+1}+\tau)}\varphi(y_{\alpha_1,1})e_{\alpha_1,1}e_{\alpha_2,2}\ldots e_{\alpha_i,i} \eqno(7)$$
Let $$f_{\alpha_1,i}=\sum\Sb1\leqslant\alpha_2\leqslant p_2\\ \\\dots\dots\dots\dots\\ \\1\leqslant\alpha_i\leqslant p_i\endSb\frac{\botshave{\prod_{1\leqslant\nu<i}}\theta(y_{\alpha_\nu,\nu}-y_{\alpha_{\nu+1},\nu+1}+v_\nu-v_i+\tau)}{\botshave{\prod_{1\leqslant\nu<i}}\theta(y_{\alpha_\nu,\nu}-y_{\alpha_{\nu+1},\nu+1}+\tau)}e_{\alpha_1,1}e_{\alpha_2,2}\ldots e_{\alpha_i,i} \eqno(8)$$
From (7) and (8) we have:
$$z_{\varphi,i}=\sum_{1\leqslant\alpha\leqslant p_1}\varphi(y_{\alpha,1})f_{\alpha,i} \eqno(9)$$
Using (4) it is easy to check the following relations between $\{f_{\alpha,i}, y_{\alpha,1}, v_i; 1\leqslant\alpha\leqslant p_1; 1\leqslant i\leqslant h\}$:
$$f_{\alpha,i}f_{\beta,j}=\frac{e^{2\pi i\tau}\theta(v_i-v_j)\theta(y_{\alpha,1}-y_{\beta,1}-2\tau)}{\theta(v_i-v_j+2\tau)\theta(y_{\alpha,1}-y_{\beta,1})}f_{\beta,j}f_{\alpha,i}+$$
$$\frac{e^{2\pi i(v_i-v_j-\tau)}\theta(2\tau)\theta(y_{\alpha,1}-y_{\beta,1}+v_j-v_i)}{\theta(v_i-v_j+2\tau)\theta(y_{\alpha,1}-y_{\beta,1})}f_{\alpha,j}f_{\beta,i}, \text{ here } \alpha\ne\beta, i<j; \eqno(10)$$
$$f_{\alpha,i}f_{\beta,i}=-e^{2\pi i(y_{\beta,1}-y_{\alpha,1})}\frac{\theta(y_{\alpha,1}-y_{\beta,1}-2\tau)}{\theta(y_{\beta,1}-y_{\alpha,1}-2\tau)}f_{\beta,i}f_{\alpha,i}, \text{ here } \alpha\ne\beta;$$
$$f_{\alpha,j}f_{\alpha,i}=e^{2\pi i\tau}f_{\alpha,i}f_{\alpha,j}, \text{ here } i<j;$$
$$f_{\alpha,i}y_{\beta,1}=y_{\beta,1}f_{\alpha,i}, \text{ here } \alpha\ne\beta; f_{\alpha,i}y_{\alpha,1}=(y_{\alpha,1}+2\tau)f_{\alpha,i};$$
$$f_{\alpha,i}v_j=(v_j-2\tau)f_{\alpha,i}, \text{ here } i\ne j; f_{\alpha,i}v_i=(v_i-4\tau)f_{\alpha,i};$$
$$[y_{\alpha,1},y_{\beta,1}]=[y_{\alpha,1},v_i]=[v_i,v_j]=0.$$
The proposition follows from the relations (10), the formula (9) and the results of [2].

\medskip

{\bf 4.3. The elliptic deformation of the algebra functions on Grassmannian.}

Let us set $\Delta=A_h$, $(\delta_i,\delta_i)=2$, $1\leqslant i\leqslant h$; $(\delta_i,\delta_{i+1})=-1$, $1\leqslant i<h$; $(\delta_i,\delta_j)=0$ if $|i-j|>1$ as above. Let $n(\delta_1)=-1$, $n(\delta_i)=0$ for $2\leqslant i\leqslant h-1$, $n(\delta_h)=m>0$. The space $Q_{\delta_1+2\delta_2+\ldots+h\delta_h}$ consists of the following functions:
$$f(x_{1,1};x_{1,2},x_{2,2};\dots;x_{1,h},\dots,x_{h,h};u_1,\dots,u_h)=$$
$$\frac{\botshave{\prod_{1\leqslant i<h}}\theta(\botshave{\sum_{1\leqslant\alpha\leqslant i}}x_{\alpha,i}-\botshave{\sum_{1\leqslant\alpha\leqslant i+1}}x_{\alpha,i+1}+\botshave{\sum_{1\leqslant\alpha\leqslant i}}u_\alpha)\botshave{\prod\Sb1\leqslant i\leqslant h,\\ 1\leqslant \alpha\ne\beta\leqslant i\endSb}\theta(x_{\alpha,i}-x_{\beta,i}-\tau)}{\botshave{\prod\Sb1\leqslant i<h,\\ 1\leqslant\alpha\leqslant i,1\leqslant\beta\leqslant i+1\endSb}\theta(x_{\alpha,i}-x_{\beta,i+1})}\times \eqno(11)$$
$$\varphi(x_{1,h},\dots,x_{h,h};u_1,\dots,u_h)$$

Here $\varphi(x_1,\dots,x_h;u_1,\dots,u_h)$ is holomorphic and symmetric in the variables $\{x_1,\dots,x_h\}$ and satisfies relations:
$$\varphi(x_1+1,x_2,\dots,x_h;u_1,\dots,u_h)=\varphi(x_1,\dots,x_h;u_1,\dots,u_h)$$
$$\varphi(x_1+\eta,x_2,\dots,x_h;u_1,\dots,u_h)=$$
$$e^{-2\pi i((m-2h+2)x_1+x_2+\ldots+x_h+u_1+\ldots+u_h)}\varphi(x_1,\dots,x_h;u_1,\dots,u_h)$$
It is easy to check that the space of such functions is isomorphic to $\Lambda^h\Theta_{m-h+1,u_1+\ldots+u_h}(\Gamma)$.

{\bf Proposition 9.} {\it For $\tau=0$ the subalgebra $\botshave{\bigoplus_{\alpha\geqslant 1}}Q_{\alpha(\delta_1+2\delta_2+\ldots+h\delta_h)}$ of the algebra $Q_{n,\Delta}(\Cal E,\tau)$ is an algebra of functions on the manifold of $h$-dimensional subspaces in $m-h+1$-dimensional vector space. This means that the algebra is commutative and generated by $Q_{\delta_1+2\delta_2+\ldots+h\delta_h}\cong \Lambda^h\Theta_{m-h+1,u_1+\ldots+u_h}(\Gamma)$ and the Plucker relations. For generic $\tau$ the algebra  $\botshave{\bigoplus_{\alpha\geqslant 1}}Q_{\alpha(\delta_1+2\delta_2+\ldots+h\delta_h)}$ is a flat deformation of the algebra of functions on the Grassmannian $Gr(h,m-h+1)$.}

{\bf Proof.} For $\tau=0$ the proof is similar to the proof of proposition 8. From our definitions it follows that $\roman{dim}_{Q_0}Q_{\alpha(\delta_1+2\delta_2+\ldots+h\delta_h)}$ does not depend on $\tau$. So this dimension is the same as for $\tau=0$.

\medskip

{\bf 4.4. Algebra $Q_{hm+1,h}(\Cal E,\tau)$.}

Let $\Delta=A_h$, $(\delta_i,\delta_i)=2$, $1\leqslant i\leqslant h$; $(\delta_i,\delta_{i+1})=-1$, $1\leqslant i<h$; $(\delta_i,\delta_j)=0$ if $|i-j|>1$ as above. Let $n(\delta_1)=1$, $n(\delta_i)=0$ for $2\leqslant i\leqslant h-1$, $n(\delta_h)=m>0.$ From the proposition 7 it follows that the algebra $Q_{n,\Delta}(\Cal E,0)$ is the algebra of polynomials in $h(m+1)$ variables: one variable of degree $\delta_1+\ldots+\delta_i$ for each $1\leqslant i<h$, $m$ variables of degree $\delta_i+\ldots+\delta_h$ for each $1<i\leqslant h$ and $m+1$ variables of degree $\delta_1+\ldots+\delta_h$. It is easy to see from this that the subalgebra $\botshave{\bigoplus_{\alpha\geqslant1}}Q_{\alpha(\delta_1+\ldots+\delta_h)}$ of the algebra $Q_{n,\Delta}(\Cal E,0)$ is the algebra of polynomials in $hm+1$ variables of degree $\delta_1+\ldots+\delta_h$. So the  subalgebra $\botshave{\bigoplus_{\alpha\geqslant1}}Q_{\alpha(\delta_1+\ldots+\delta_h)}$ of the algebra $Q_{n,\Delta}(\Cal E,\tau)$ is a flat deformation of the algebra of polynomials in $hm+1$ variables. The Hilbert function is $1+\botshave{\sum_{\alpha\geqslant1}}\roman{dim}_{Q_0}Q_{\alpha(\delta_1+\ldots+\delta_h)} w^\alpha=(1-w)^{-(hm+1)}$.

{\bf Proposition 10.} {\it The subalgebra $\botshave{\bigoplus_{\alpha\geqslant1}}Q_{\alpha(\delta_1+\ldots+\delta_h)}$ of the algebra $Q_{n,\Delta}(\Cal E,\tau)$ is isomorphic to $Q_{hm+1,h}(\Cal E,\frac{2}{hm+1}\tau)$.} 

{\bf Proof.} The proof is similar to the proof of proposition 8, using the construction of the algebras $Q_{n,k}(\Cal E,\tau)$ (see [2]).

\medskip

{\bf 4.5. Generalized elliptic $R$-matrix.}

Let $\Delta=A_h$, $(\delta_i,\delta_i)=2$, $1\leqslant i\leqslant h$; $(\delta_i,\delta_{i+1})=-1$, $1\leqslant i<h$; $(\delta_i,\delta_j)=0$ if $|i-j|>1$ as above. We fix $\nu$, $1<\nu<h$. Let $n(\delta_\nu)=1$, $n(\delta_i)=0$ for $1\leqslant i\leqslant h$, $i\ne\nu$. According to the proposition 7 the algebra $Q_{n,\Delta}(\Cal E,0)$ is the algebra of polynomials in $\nu(h-\nu+1)$ variables: one variable of degree $\delta_{i_1}+\delta_{i_1+1}+\ldots+\delta_{i_2}$ for each $i_1$, $i_2$ such that $1\leqslant i_1\leqslant\nu\leqslant i_2\leqslant h$. The algebra $Q_{n,\Delta}(\Cal E,\tau)$ is a flat deformation of this polynomial ring. It is clear that in the algebra $Q_{n,\Delta}(\Cal E,\tau)$ we have $\roman{dim}_{Q_0}Q_{\delta_{i_1}+\ldots+\delta_{i_2}}=1$ for $1\leqslant i_1\leqslant\nu\leqslant i_2\leqslant h$, $Q_{\delta_{i_1}+\ldots+\delta_{i_2}}=0$ for $\nu<i_1$ or $i_2<\nu$. Let $e_{i_1,i_2}$ be a nonzero element from $Q_{\delta_{i_1}+\delta_{i_1+1}+\ldots+\delta_{i_2}}$ (for $1\leqslant i_1\leqslant\nu\leqslant i_2\leqslant h$). It is easy to check that we can choose $e_{i_1,i_2}$ in such a form:
$$e_{i_1,i_2}=$$
$$\frac{\botshave{\prod_{i_1\leqslant\varphi<\nu}}\theta(x_{1,\varphi}-x_{1,\varphi+1}+\botshave{\sum_{i_1\leqslant\alpha\leqslant\varphi}}u_{\alpha}-\tau)\botshave{\prod_{\nu\leqslant\psi<i_2}}\theta(x_{1,\psi}-x_{1,\psi+1}-\botshave{\sum_{\psi<\alpha\leqslant i_2}}u_{\alpha}+\tau)}{\botshave{\prod_{i_1\leqslant p<i_2}}\theta(x_{1,p}-x_{1,p+1})}\times \eqno(12)$$
$\theta(x_{1,\nu}+\botshave{\sum_{i_1\leqslant\alpha\leqslant i_2}}u_{\alpha}-2\tau)$

Let us denote $v_i=u_i+u_{i+1}+\ldots+u_\nu$ for $1\leqslant i\leqslant\nu$ and $w_i=u_\nu+u_{\nu+1}+\ldots+u_i$ for $\nu\leqslant i\leqslant h$. We denote $q(j,j^\prime)=e^{2\pi i\tau}$ for $j<j^\prime$, $q(j,j^\prime)=e^{-2\pi i\tau}$ for $j>j^\prime$, $q(j,j)=1$.

{\bf Proposition 11.} {\it The elements $\{e_{i_1,i_2}, v_{i_1}, w_{i_2}; 1\leqslant i_1\leqslant\nu\leqslant i_2\leqslant h\}$ satisfy following commutation relations:
$$e_{i_1,i_2}e_{i_1^\prime,i_2^\prime}=\frac{\theta(v_{i_1}-v_{i_1^\prime}-2\tau)\theta(w_{i_2^\prime}-w_{i_2})q(i_1^\prime,i_1)q(i_2^\prime,i_2)}{\theta(v_{i_1}-v_{i_1^\prime})\theta(w_{i_2^\prime}-w_{i_2}-2\tau)}e_{i_1^\prime,i_2^\prime}e_{i_1,i_2}+ \eqno(13)$$
$$\frac{\theta(-2\tau)q(i_2^\prime,i_2)\theta(w_{i_2^\prime}-w_{i_2}-v_{i_1}+v_{i_1^\prime})}{\theta(v_{i_1^\prime}-v_{i_1})\theta(w_{i_2^\prime}-w_{i_2}-2\tau)}e_{i_1,i_2^\prime}e_{i_1^\prime,i_2},\text{ for } i_1\ne i_1^\prime, i_2\ne i_2^\prime;$$
$$e_{i_1,i_2}e_{i_1^\prime,i_2}=q(i_1^\prime,i_1)e_{i_1^\prime,i_2}e_{i_1,i_2}; e_{i_1,i_2}e_{i_1,i_2^\prime}=q(i_2^\prime,i_2)e_{i_1,i_2^\prime}e_{i_1,i_2};$$
$$e_{i_1,i_2}v_{i_1^\prime}=v_{i_1^\prime}e_{i_1,i_2}\text{ for }i_1\ne i_1^\prime; e_{i_1,i_2}w_{i_2^\prime}=w_{i_2^\prime}e_{i_1,i_2}\text{ for }i_2\ne i_2^\prime;$$
$$e_{i_1,i_2}v_{i_1}=(v_{i_1}-2\tau)e_{i_1,i_2}; e_{i_1,i_2}w_{i_2}=(w_{i_2}-2\tau)e_{i_1,i_2};$$
$$[v_{i_1},v_{i_1^\prime}]=[w_{i_2},w_{i_2^\prime}]=[v_{i_1},w_{i_2}]=0$$}

{\bf Proof.} It is clear that $$e_{i_1,i_2}e_{i_1^\prime,i_2^\prime},\; e_{i_1^\prime,i_2^\prime}e_{i_1,i_2},\; e_{i_1,i_2^\prime}e_{i_1^\prime,i_2}\in Q_{\delta_{i_1}+\delta_{i_1+1}+\ldots+\delta_{i_2}+\delta_{i_1^\prime}+\delta_{i_1^\prime+1}+\ldots+\delta_{i_2^\prime}}$$ but $\roman{dim}_{Q_0}Q_{\delta_{i_1}+\ldots+\delta_{i_2}+\delta_{i_1^\prime}+\ldots+\delta_{i_2^\prime}}=2$, so there is a linear relation between these elements. It is easy to find the coefficients of this relation using (2). Another relations in (13) are easy to check directly.

{\bf Remark.} If we put $n(\delta_\nu)=m>0$ instead of $n(\delta_\nu)=1$, then we will have an algebra that is a deformation of the polynomial ring in $m\nu(h-\nu+1)$ variables: $m$ variables of degree $\delta_{i_1}+\delta_{i_1+1}+\ldots+\delta_{i_2}$ for each $1\leqslant i_1\leqslant\nu\leqslant i_2\leqslant h$. It is clear that $\roman{dim}_{Q_0}Q_{\delta_{i_1}+\ldots+\delta_{i_2}}=m$ for $1\leqslant i_1\leqslant\nu\leqslant i_2\leqslant h$, $Q_{\delta_{i_1}+\ldots+\delta_{i_2}}=0$ for $i_1>\nu$ or $i_2<\nu$ and $\roman{dim}_{Q_0}Q_{\delta_{i_1}+\ldots+\delta_{i_2}+\delta_{i_1^\prime}+\ldots+\delta_{i_2^\prime}}=2m^2$ for $1\leqslant i_1,i_1^\prime\leqslant\nu\leqslant i_2,i_2^\prime\leqslant h$. So the elements from the product $Q_{\delta_{i_1}+\ldots+\delta_{i_2}}*Q_{\delta_{i_1^\prime}+\ldots+\delta_{i_2^\prime}}$ are uniquely expressed as linear combination of the elements from the space $Q_{\delta_{i_1^\prime}+\ldots+\delta_{i_2^\prime}}*Q_{\delta_{i_1}+\ldots+\delta_{i_2}}\oplus Q_{\delta_{i_1}+\ldots+\delta_{i_2^\prime}}*Q_{\delta_{i_1^\prime}+\ldots+\delta_{i_2}}$. So we have a map $$Q_{\delta_{i_1}+\ldots+\delta_{i_2}}\otimes Q_{\delta_{i_1^\prime}+\ldots+\delta_{i_2^\prime}}\to Q_{\delta_{i_1^\prime}+\ldots+\delta_{i_2^\prime}}\otimes Q_{\delta_{i_1}+\ldots+\delta_{i_2}}\oplus Q_{\delta_{i_1}+\ldots+\delta_{i_2^\prime}}\otimes Q_{\delta_{i_1^\prime}+\ldots+\delta_{i_2}}$$
 It is clear that this map satisfies the Yang-Baxter equations.

\newpage

\centerline{\bf \S5. Case of affine root system.}
\centerline{\bf A family of commuting elements.}
\medskip

Let $\Delta$ be the root system of the Lie algebra $\widehat{sl_h}$. So $(\delta_i,\delta_i)=2$, $(\delta_i,\delta_{i+1})=-1$ and $(\delta_i,\delta_j)=0$ otherwise. Here we assume that $i\in \Bbb Z_h$. It is possible to prove that for dominant $n$ the Hilbert function of $Q_{n,\Delta}(\Cal E,\tau)$ is equal to $$\botshave{\prod_{\gamma\in \Delta^+}}(1-w^\gamma)^{-n(\gamma)k(\gamma)}\botshave{\prod_{\alpha\geqslant\mu}}(1-w^{\alpha(\delta_1+\ldots+\delta_h)})^{-\alpha n(\delta_1+\ldots+\delta_h)}$$ Here $k(\gamma)$ is a multiplicity of the root $\gamma$ in the algebra $\widehat{sl_h}$, $\mu=2$ for $h=2$ and $\mu=1$ for $h>2$. So the algebra $Q_{n,\Delta}(\Cal E,\tau)$ is bigger then in the finite dimensional case. Note that the element $u_1+\ldots+u_h$ is central in the algebra $Q_{n,\Delta}(\Cal E,\tau)$, because the root $\delta_1+\ldots+\delta_h$ is imaginary. We denote by $Q_{n,\Delta,c}(\Cal E,\tau)$ (for $c\in\Bbb C$) the factor algebra of $Q_{n,\Delta}(\Cal E,\tau)$ by the relation $u_1+\ldots+u_h=c$. If $n$ is dominant and $n\ne0$ then the Hilbert function of the algebra $Q_{n,\Delta,c}(\Cal E,\tau)$ is the same as for the algebra $Q_{n,\Delta}(\Cal E,\tau)$ and does not depend on $c$. Let us $n=0$. The Hilbert function of the algebra $Q_{0,\Delta,c}(\Cal E,\tau)$ is equal to 1 if $c\notin\Gamma$ and is equal to $\botshave{\prod_{\alpha\geqslant1}}(1-w^{\alpha(\delta_1+\ldots+\delta_h)})^{-h}$ for $h>2$ and $(1-w_1w_2)^{-1}\botshave{\prod_{\alpha\geqslant2}}(1-w_1^\alpha w_2^\alpha)^{-2}$ for $h=2$ if $c\in\Gamma$.

{\bf Proposition 12.} {\it Let $c\in\Gamma$. The algebra $Q_{0,\Delta,c}(\Cal E,\tau)$ is commutative for each $\tau\in\Cal E$ and isomorphic to polynomials in infinite number of variables. Namely, for $h>2$ it is the algebra of polynomials in $h$ variables of degree $\alpha(\delta_1+\ldots+\delta_h)$ for each $\alpha\geqslant1$ and for $h=2$ it is the algebra of polynomials in one variable of degree $\delta_1+\delta_2$ and two variables of degree $\alpha(\delta_1+\delta_2)$ for each $\alpha\geqslant2$.}

{\bf Proof.} Let $c=0$. Firstly we consider the case $h>2$. Let $g(z_1,\dots,z_h)$ be the holomorphic function in variables $\{z_j; j\in\Bbb Z_h\}$ that satisfies the relations:

 $g(z_1,\dots,z_j+1,\dots,z_h)=g(z_1,\dots,z_h),$

 $g(z_1,\dots,z_j+\eta,\dots,z_h)=e^{-2\pi i(2z_j-z_{j-1}-z_{j+1}-\eta+u_j)}g(z_1,\dots,z_h)$ for each $j\in\Bbb Z_h$. It is easy to see that the space of such functions is $h$-dimensional (over field of functions in variables $\{u_j\}$) and $g(z_1+p,\dots,z_h+p)=g(z_1,\dots,z_h)$ for each $p\in\Bbb C$. Let $K_{\alpha,g}$ be the following element of the space $Q_{\alpha(\delta_1+\ldots+\delta_h)}$:
$$K_{\alpha,g}=\frac{\botshave{\prod\Sb i\in\Bbb Z_h,\\ 1\leqslant\mu,\nu\leqslant\alpha\endSb}\theta(x_{\mu,i}-x_{\nu,i}-2\tau)\cdot g(x_{1,1}+\ldots+x_{\alpha,1},\dots,x_{1,h}+\ldots+x_{\alpha,h})}{\botshave{\prod\Sb i\in\Bbb Z_h,\\ 1\leqslant\mu,\nu\leqslant\alpha\endSb}\theta(x_{\mu,i}-x_{\nu,i+1})} \eqno(14)$$

One can check by direct calculation that $K_{\alpha_1,g_1}*K_{\alpha_2,g_2}=K_{\alpha_2,g_2}*K_{\alpha_1,g_1}$ for each $\alpha_1,\alpha_2$; $g_1,g_2$.

For $h=2$ the elements $K_{\alpha,g}$ defined by (14) commute also, but in this case $K_{\alpha,g}\notin Q_{\alpha(\delta_1+\delta_2)}$ because $K_{\alpha,g}$ has a pole of order 2. In this case the algebra $Q_{0,\Delta,0}(\Cal E,\tau)$ is a subalgebra of the algebra generated by $\{K_{\alpha,g}\}$.

\newpage

\centerline{\bf \S6. Elliptic deformation of the Poison algebra.}
\centerline{\bf Another family of commuting elements.}
\medskip

Let $\tau_1,\tau_2,\tau_3\in\Cal E$, $\tau_1+\tau_2+\tau_3=0$ and $m\in\Bbb N$. We define an associative $\Bbb N$-graded algebra $H_m(\Cal E,\tau_1,\tau_2,\tau_3)$ by the following: as a linear space $H_m(\Cal E,\tau_1,\tau_2,\tau_3)=F_1\oplus F_2\oplus F_3\oplus\ldots$, where $F_\alpha*F_\beta\subset F_{\alpha+\beta}$ and $F_\alpha$ is a space of functions $f(x_1,\dots,x_\alpha)$ on $\Bbb C^\alpha$ that satisfies the properties:

1. $f(x_1,\dots,x_\alpha)$ is symmetric by $x_1,\dots,x_\alpha$.

2. $f(x_1,\dots,x_\alpha)$ is holomorphic outside the divisors $\{x_\mu-x_\nu=0\}$ and has a pole of order $\leqslant2$ on these divisors.

3. For $\alpha\geqslant3$ we have:$$f(x,x+\tau_1,x+\tau_1+\tau_2,x_4,\dots,x_\alpha)=0$$
$$f(x,x+\tau_2,x+\tau_2+\tau_1,x_4,\dots,x_\alpha)=0$$

4. $f(x_1+1,x_2,\dots,x_\alpha)=f(x_1,\dots,x_\alpha)$,

$f(x_1+\eta,x_2,\dots,x_\alpha)=e^{-2\pi i(mx_1+c)}f(x_1,\dots,x_\alpha), \text{ here }c\in\Bbb C \text{ is fixed. }$

We define the product $*$ in the algebra $H_m(\Cal E,\tau_1,\tau_2,\tau_3)$ by the following rule: for $f\in F_\alpha$, $g\in F_\beta$ we set:
$$f*g(x_1,\dots,x_{\alpha+\beta})=$$
$$\frac{1}{\alpha!\beta!}\sum_{\sigma\in S_{\alpha+\beta}}f(x_{\sigma_1},\dots,x_{\sigma_\alpha})g(x_{\sigma_{\alpha+1}},\dots,x_{\sigma_{\alpha+\beta}})\prod\Sb1\leqslant\mu\leqslant\alpha,\\ \alpha+1\leqslant\nu\leqslant\alpha+\beta\endSb\lambda(x_{\sigma_\mu},x_{\sigma_\nu}),$$
here $\lambda(x,y)=\frac{\theta(x-y-\tau_1)\theta(x-y-\tau_2)\theta(x-y-\tau_3)}{\theta(x-y)^3}$.

It is clear that for $m>0$ the algebra $H_m(\Cal E,\tau_1,\tau_2,\tau_3)$ does not depend on $c$. It is possible to proof that for $m>0$ the algebra $H_m(\Cal E,\tau_1,\tau_2,\tau_3)$ is a deformation in class of $\Bbb N$-graded associative algebras of the universal enveloping algebra of the following $\Bbb N$-graded Lie algebra $h_{m,c}(\Cal E)$: as a vector space, $h_{m,c}(\Cal E)=\Theta_{m,c}(\Gamma)\oplus\Theta_{2m,2c}(\Gamma)\oplus\ldots$ and for $\varphi\in\Theta_{\alpha m,\alpha c}(\Gamma)$, $\psi\in\Theta_{\beta m,\beta c}(\Gamma)$, the commutator $[\varphi,\psi]=\beta\varphi^\prime\psi-\alpha\psi^\prime\varphi\in\Theta_{(\alpha+\beta)m,(\alpha+\beta)c}(\Gamma)$. The Hilbert function of the algebra $H_m(\Cal E,\tau_1,\tau_2,\tau_3)$ is $1+\botshave{\sum_{\alpha\geqslant1}}\roman{dim}\,F_\alpha w^\alpha=\botshave{\prod_{\alpha\geqslant1}}(1-w^\alpha)^{-m\alpha}$ (see also [5]).

Now let us apply this construction to the case $m=0$. For $c\notin\Gamma$ we have $H_0(\Cal E,\tau_1,\tau_2,\tau_3)=0$. For $c\in\Gamma$ the Hilbert function of the algebra $H_0(\Cal E,\tau_1,\tau_2,\tau_3)$ is equal to $\botshave{\prod_{\alpha\geqslant1}}(1-w^\alpha)^{-1}$. In this case the algebra  $H_0(\Cal E,\tau_1,\tau_2,\tau_3)$ is commutative for each $\tau_1,\tau_2,\tau_3\in\Bbb C$, $\tau_1+\tau_2+\tau_3=0$ and isomorphic to the polynomial algebra in infinite number of variables: one variable of degree $\alpha$ for each $\alpha\geqslant1$. Let us $c=0$. We define the elements $K_\alpha\in F_\alpha$ by the formula:
$$K_\alpha(x_1,\dots,x_\alpha)=\prod_{1\leqslant\mu<\nu\leqslant\alpha}\frac{\theta(x_\mu-x_\nu-\tau_1)\theta(x_\mu-x_\nu+\tau_1)}{\theta(x_\mu-x_\nu)^2} \text{ for }\alpha>1,$$
$$K_1(x_1)=1$$

It is easy to check by direct calculation that $K_\alpha*K_\beta=K_\beta*K_\alpha$ for each $\alpha$, $\beta\in\Bbb N$ (see for details [5]).

\newpage

\centerline{\bf References.}
\medskip
1. A.V.Odesskii and B.L.Feigin, "Sklyanin's elliptic algebras", Funkts. Anal. Philozhen., 23, No. 3, 45-54 (1989).

2. A.V.Odesskii and B.L.Feigin, "Constructions of Sklyanin elliptic algebras and quantum $R$-matrices", Funkts. Anal. Philozhen. 27, No. 1, 37-45 (1993).

3. B.L.Feigin and A.V.Odesskii, Vector Bundles on Elliptic Curve and Sklyanin Algebras, RIMS-1032, September 1995, Kyoto University, Kyoto, Japan.

4. Cherednik I.V. On $R$-matrix quantization of formal loop groups. Group theoretical methods in physics, Voll. 2 (Yurmala, 1985), 161-180, VNU Sci. Press, Utrecht, 1986.

5. B.L.Feigin and A.V.Odesskii, A family of elliptic algebras. Internat. Math. Res. Notices 1997, No 11, 531-539.

6. A.V.Odesskii and B.L.Feigin, Elliptic deformation of current algebras and their representations by difference operators. Funct. Anal. Appl. 31 (1997), No 3, 193-203 (1998).

\enddocument